\documentclass[11pt,british]{amsart}
\usepackage{ifpdf}
\ifpdf\pdfoutput=1\pdftrue\fi\newif\pdf
\ifpdf\relax\else
\usepackage[T1]{fontenc}
\fi

\usepackage{babel,eucal,url,amssymb,stmaryrd,
enumerate,amscd,
}

\usepackage[pagebackref]{hyperref}

\theoremstyle{plain}
\newtheorem{thm}{Theorem}[section]
\newtheorem{lem}[thm]{Lemma}
\newtheorem{pro}[thm]{Proposition}
\newtheorem{co}[thm]{Corollary}

\theoremstyle{definition}
\newtheorem{defn}[thm]{Definition}

\theoremstyle{remark}
\newtheorem{rem}[thm]{Remark}

\newcommand{\Gtwo}{\ifmmode{{\rm G}_2}\else{${\rm G}_2$}\fi}

\date{\today}

\begin{document}

\title{Canonical connections on paracontact manifolds}

\author{Simeon Zamkovoy}
\address[Zamkovoy]{University of Sofia "St. Kl. Ohridski"\\
Faculty of Mathematics and Informatics,\\
Blvd. James Bourchier 5,\\
1164 Sofia, Bulgaria} \email{zamkovoy@fmi.uni-sofia.bg}

\begin{abstract}
{ The canonical paracontact connection is defined and it is shown that its torsion is the obstruction
the paracontact manifold to be paraSasakian. A $\mathcal{D}$-homothetic transformation is determined
as a special gauge transformation. The $\eta$-Einstein manifold are defined, it is
 prove that their scalar curvature is a constant and it is shown that in the paraSasakian case these
 spaces  can be obtained from  Einstein paraSasakian manifolds with a  $\mathcal{D}$-homothetic
 transformations. It is shown that an almost paracontact
structure admits a connection with totally skew-symmetric torsion
if and only if the Nijenhuis tensor of the paracontact structure
is skew-symmetric and the defining vector field is Killing.}


MSC: 53C15, 5350, 53C25, 53C26, 53B30
\end{abstract}

\maketitle \setcounter{tocdepth}{3} \tableofcontents

\section{Introduction}

In \cite{K1} Kaneyuki and Konzai defined the almost paracontact
structure on pseudo-Riemannian manifold $M$ of dimension $(2n+1)$
and constructed the almost paracomplex structure on
$M^{(2n+1)}\times \mathbb{R}$.
In this paper we study the properties of an almost paracontact
metric manifold.
We consider gauge (conformal) transformations of a paracontact
manifold i.e. transformations preserving the paracontact
structure.  We define $\mathcal{D}$-homothetic transformations as
a special gauge transformation (homothetic) and study the behavior
of the Einstein condition under  $\mathcal{D}$-homothetic
transformations on a paracontact metric manifold. We consider the
$\eta$-Einstein manifold,
 prove that their scalar curvature is a constant and show that in the paraSasakian case these spaces
 are the images of  Einstein paraSasakian manifolds under  $\mathcal{D}$-homothetic transformations.

We define a canonical paracontact connection on a paracontact metric manifold
which seems to be the paracontact analogue of the (generalized) Tanaka-Webster connection.
We show that the torsion of this connection vanishes exactly when
the structure is para-Sasakian and compute the gauge
transformation of its scalar curvature.

We introduce and study also the notion of paracontact manifolds
with torsion. The paracontact manifolds with torsion are
manifolds, which admit a linear almost paracontact connection with
totaly skew-symmetric torsion. We prove that an almost paracontact
structure admits a connection with totally skew-symmetric torsion
if and only if the Nijenhuis tensor of the paracontact structure
is skew-symmetric and the defining vector field is Killing. In the
contact case this connection is studied in \cite{F1}.

\section{Almost paracontact manifolds}
A (2n+1)-dimensional smooth manifold $M^{(2n+1)}$
has an \emph{almost paracontact structure} 
$(\varphi,\xi,\eta)$
if it admits a tensor field
$\varphi$ of type $(1,1)$, a vector field $\xi$ and a 1-form
$\eta$ satisfying the  following compatibility conditions 
\begin{eqnarray}
  \label{f82}
    & &
    \begin{array}{cl}
          (i)   & \varphi(\xi)=0,\quad \eta \circ \varphi=0,\quad
          \\[5pt]
          (ii)  & \eta (\xi)=1 \quad \varphi^2 = id - \eta \otimes \xi,
          \\[5pt]
          (iii) & \textrm{let $\mathbb D=Ker~\eta$ be the  horizontal distribution generated
                          by $\eta$, then}
          \\[3pt]
                & \textrm{the tensor field $\varphi$ induces an almost paracomplex structure (see \cite{K2})}
          \\[3pt]
                &  \textrm{on each fibre on $\mathbb D$.}
    \end{array}
\end{eqnarray}
Recall that an almost paracomplex structure on an 2n-dimensional
manifold is a (1,1)-tensor $J$ such that $J^2=1$ and the
eigensubbundles $T^+,T^-$ corresponding to the eigenvalues $1,-1$
of $J$, respectively have equal dimension $n$. The Nijenhuis
tensor $N$ of $J$, given by
$N_{J}(X,Y)=[JX,JY]-J[JX,Y]-J[X,JY]+[X,Y],$ is the obstruction for
the integrability of the eigensubbundles $T^+,T^-$. If $N=0$ then
the almost paracomplex structure is called paracomplex or
integrable.

An immediate consequence of the definition of the almost
paracontact structure is that the endomorphism $\varphi$ has rank
$2n$,
$\varphi \xi=0$ and $\eta \circ \varphi=0$, (see \cite{B1,B2}
for the almost contact case).


If a manifold $M^{(2n+1)}$ with $(\varphi,\xi,\eta)$-structure
admits a pseudo-Riemannian metric $g$ such that
\begin{equation}\label{con}
g(\varphi X,\varphi Y)=-g(X,Y)+\eta (X)\eta (Y),
\end{equation}
then we say that $M^{(2n+1)}$ has an almost paracontact metric structure and
$g$ is called \emph{compatible} metric. Any compatible metric $g$ with a given almost paracontact
structure is necessarily of signature $(n+1,n)$.

Setting $Y=\xi$, we have
$\eta(X)=g(X,\xi).$

Any almost paracontact structure admits a compatible metric.
Indeed, if G is any metric, first set
$\overline{G}(X,Y)=G(\varphi^2 X,\varphi^2 Y)+\eta(X)\eta(Y)$;
then $\eta(X)=\overline{G}(X,\xi)$. Now define $g$ by
$g(X,Y)=\frac{1}{2}(\overline{G}(X,Y)-\overline{G}(\varphi X,\varphi Y)+\eta(X)\eta(Y))$
end check $g$ is compatible.

The fundamental 2-form
\begin{equation}\label{fund}
F(X,Y)=g(X,\varphi Y)
\end{equation}
is non-degenerate on the horizontal distribution $\mathbb D$ and
$\eta\wedge F^n\not=0$.
\begin{defn}
If $g(X,\varphi Y)=d\eta(X,Y)$ (where
$d\eta(X,Y)=\frac12(X\eta(Y)-Y\eta(X)-\eta([X,Y])$ then $\eta$ is
a paracontact form and the almost paracontact metric manifold
$(M,\varphi,\eta,g)$ is said to be $\emph{paracontact metric
manifold}$.
\end{defn}


The manifold $M$ is orientable exactly when the canonical line bundle $E=\{\eta\in \Lambda^1: Ker~\eta=\mathbb D\}$
is orientable, since $\mathbb D$ is orientable by the paracomplex structure $\varphi$. Any two contact forms
$\bar\eta,\eta\in E$ are connected by
\begin{equation}\label{gauge}
\bar\eta=\sigma\eta,
\end{equation}
where $\sigma$ is non-vanishing smooth function on $M$. We study
this conformal (gauge) transformation in $Section~\ref{gau}$

\begin{rem}
We mention that some authors say $M^{(2n+1)}$ has an almost
paracontact metric structure if it admits a Riemannian metric $g$
such that $g(\varphi X,\varphi Y)=g(X,Y)-\eta (X)\eta (Y)$ (see
\cite{R1,S1}). In our paper the metric is a pseudo-Riemannian and
metric satisfies a condition $(\ref{con})$
\end{rem}

For a manifold $M^{(2n+1)}$ with an almost paracontact metric
structure $(\varphi,\xi,\eta,g)$ we can also construct a useful
local orthonormal basis. Let U be a coordinate neighborhood on M
and $X_1$ any unit vector field on U orthogonal to $\xi$. Then
$\varphi X_1$ is a vector field orthogonal to  both X and $\xi$,
and $|\varphi X_1|^2=-1$. Now choose a unit vector field $X_2$
orthogonal to $\xi$, $X_1$ and $\varphi X_1$. Then $\varphi X_2$ is
also vector field orthogonal to $\xi$, $X_1$, $\varphi X_1$ and
$X_2$, and $|\varphi X_2|^2=-1$. Proceeding in this way we obtain
a local orthonormal basis $(X_i,\varphi X_i,\xi),i=1...n$ called a
\emph{$\varphi$-basis}.

Hence, an almost paracontact metric manifold
$(M^{2n+1},\varphi,\eta,\xi,g)$ is an odd dimensional manifold
with a structure group $\mathbb U(n,\mathbb R)\times Id,$, where
$\mathbb U(n,\mathbb R)$ is the para-unitary group isomorphic to
$\mathbb {GL}(n,\mathbb R)$.



Let $M^{(2n+1)}$ be an almost paracontact manifold with structure
$(\varphi,\xi,\eta)$ and consider the manifold $M^{(2n+1)}\times
\mathbb R$. We denote a vector field on $M^{(2n+1)}\times \mathbb
R$ by $(X,f\frac{d}{dt})$ where X is tangent to $M^{(2n+1)}$, t is
the coordinate on $\mathbb R$ and f is a $C^{\infty}$ function on
$M^{(2n+1)}\times \mathbb R$. An almost
paracomplex structure J on $M^{(2n+1)}\times \mathbb R$ is defined in \cite{K1}
by $$J(X,f\frac{d}{dt})=(\varphi X+f\xi,\eta (X)\frac{d}{dt}).$$
If J is integrable, we say that the almost paracontact structure
$(\varphi,\xi,\eta)$ is \emph{normal}.

As the vanishing of the Nijenhuis tensor of J is necessary and
sufficient condition for integrability, we  express the
condition of normality in terms of Nijenhuis tensor of $\varphi$.
Since $N_{J}$ is tensor field of type $(1,2)$, it suffices to
compute $N_{J}((X,0),(Y,0))$ and $N_{J}((X,0),(0,\frac{d}{dt}))$
for vector fields X and Y on $M^{(2n+1)}$.

$$N_{J}((X,0),(Y,0))=([X,Y],0)+([\varphi X,\varphi Y],(\varphi X\eta (Y)-\varphi Y\eta (X))\frac{d}{dt})-$$
$$-(\varphi[\varphi X,Y]-Y\eta (X)\xi,\eta ([\varphi X,Y])\frac{d}{dt})-(\varphi[X,\varphi Y]+X\eta (Y)\xi,\eta ([X,\varphi Y])\frac{d}{dt})=$$
$$(N_{\varphi}(X,Y)-2d\eta (X,Y)\xi,((\pounds_{\varphi X}\eta)Y-(\pounds_{\varphi Y}\eta)X)\frac{d}{dt}).$$

$$N_{J}((X,0),(0,\frac{d}{dt}))=[(\varphi X,\eta (X)\frac{d}{dt}),(\xi,0)]-J[(X,0),(\xi,0)]=$$
$$=([\varphi X,\xi],-\xi(\eta(X))\frac{d}{dt})-(\varphi [X,\xi],\eta ([X,\xi])\frac{d}{dt})=$$
$$=-((\pounds_{\xi}\varphi)X,(\pounds_{\xi}\eta)X\frac{d}{dt}).$$

We are thus lead to define tensors $N^{(1)}$, $N^{(2)}$, $N^{(3)}$
and $N^{(4)}$ by
$$N^{(1)}(X,Y)=N_{\varphi}(X,Y) -2d\eta (X,Y)\xi,$$
$$N^{(2)}(X,Y)=(\pounds_{\varphi X}\eta)Y-(\pounds_{\varphi Y}\eta)X,$$
$$N^{(3)}(X)=(\pounds_{\xi}\varphi)X,$$
$$N^{(4)}(X)=(\pounds_{\xi}\eta)X.$$
Clearly the almost paracontact structure $(\varphi,\xi,\eta)$ is
normal if and only if these four tensors vanish.
\begin{pro}\label{t2}
For an almost paracontact structure $(\varphi,\xi,\eta)$ the
vanishing of $N^{(1)}$ implies the vanishing $N^{(2)}$, $N^{(3)}$
and $N^{(4)}$;

For a paracontact structure $(\varphi,\xi,\eta,g)$, $N^{(2)}$ and
$N^{(4)}$ vanish. Moreover $N^{(3)}$ vanishes if and only if $\xi$
is a Killing vector field.
\end{pro}
\begin{proof}
Setting $Y=\xi$ in
$d\eta(X,Y)=\frac12(X\eta(Y)-Y\eta(X)-\eta([X,Y])$ and we get
$d\eta(X,\xi)=0$. We have
$$0=N_{\varphi}(X,\xi)=-\varphi[\varphi X,\xi]+\varphi^2[X,\xi]=\varphi((\pounds_{\xi}\varphi)X).$$
Applying $\varphi$ and noting that $d\eta(\varphi X,\xi)=0$
implies $\eta([\xi,\varphi X])=0$, we have $N^{(3)}=0$. Moreover
$(\pounds_{\xi}\eta)\varphi X=0$, but $(\pounds_{\xi}\eta)\xi=0$
is immediate and hance $N^{(4)}=0$. Finally, we have
$$N_{\varphi}(\varphi X,Y)-2d\eta(\varphi X,Y)\xi=-N^{(2)}(X,Y)\xi$$
which simplifies to $N^{(2)}=0$.

If the structure is paracontact we have already seen that
$N^{(4)}(X)=(\pounds_{\xi}\eta)X=
2d\eta(\xi,X)=0.$
Now $N^{(2)}$ can be written
$$N^{(2)}(X,Y)=2d\eta(\varphi X,Y)+2d\eta(X,\varphi Y)=2g(\varphi X,\varphi Y)-2g(\varphi Y,\varphi X)=0.$$
Turning to $N^{(3)}$, since $d\eta$ is invariant under the action
of $\xi$, we have
$$0=(\pounds_{\xi}d\eta)(X,Y)=\xi g(X,\varphi Y)-g([\xi,X],\varphi Y)-g(X,\varphi [\xi,Y])=$$
$$=(\pounds_{\xi}g)(X,Y)+g(X,N^{(3)}(Y)).$$
\end{proof}
A paracontact structure for which $\xi$ is Killing vector field is
called a \emph{K-paracontact structure}.
\begin{pro}\label{l1}
For an almost paracontact metric structure $(\varphi,\xi,\eta,g)$,
the covariant derivative $\nabla\varphi$ of $\varphi$ with respect
to the Levi-Civita connection $\nabla$ is given by
\begin{gather}\label{f1}
2g((\nabla_X \varphi)Y,Z)= -dF(X,Y,Z)-dF(X,\varphi Y,\varphi
Z)-N^{(1)}(Y,Z,\varphi X)\\\nonumber
+N^{(2)}(Y,Z)\eta (X)-2d\eta (\varphi
Z,X)\eta (Y)+2d\eta (\varphi Y,X)\eta (Z).
\end{gather}
For a paracontact metric structure $(\varphi,\xi,\eta,g)$, the formula \eqref{f1}
simplifies to
\begin{equation}\label{f2}
2g((\nabla_X \varphi)Y,Z)= -N^{(1)}(Y,Z,\varphi X)-2d\eta (\varphi
Z,X)\eta (Y)+2d\eta (\varphi Y,X)\eta (Z)
\end{equation}
\end{pro}
\begin{proof}
The Levi-Civita connection $\nabla$ with respect to $g$ is given
by
$$2g(\nabla_XY,Z)=Xg(Y,Z)+Yg(Z,X)-Zg(X,Y)+g([X,Y],Z)+$$
$$+g([Z,X],Y)-g([Y,Z],X).$$
On the other hand, $dF$ can be expressed in the following way
$$dF(X,Y,Z)=XF(Y,Z)+YF(Z,X)+ZF(X,Y)-F([X,Y],Z)-$$
$$-F([Z,X],Y)-F([Y,Z],X).$$
The last two equations imply \eqref{f1}. The equation \eqref{f2}
follows from equation $(\ref{f1})$ and equalities $N^{(2)}=0$ and
$F=d\eta$.
\end{proof}
We have seen that on a contact manifold, $N^{(3)}$ vanishes if and
only if $\xi$ is Killing ($Proposition~\ref{t2}$) 
For a general paracontact
structure the tensor field $N^{(3)}$ encodes many important
properties and for simplicity we define a tensor field $h$ on a
paracontact manifold by
$$h=\frac{1}{2}\pounds_{\xi}\varphi = \frac{1}{2}N^{(3)}.$$

\begin{lem}\label{l2}
On a paracontact matric manifold, h is a symmetric operator,
\begin{equation}\label{f3}
\nabla_X\xi=- \varphi X+\varphi hX,
\end{equation}
h anti-commutes with $\varphi$ and $trh=h\xi =0$.
\end{lem}
\begin{proof}
We have already seen that on a paracontact metric manifold,
$\nabla_{\xi}\varphi =0$, $\nabla_{\xi}\xi=0$ and $N^{(2)}=0$.
Thus
$$-g((\pounds_{\xi}\varphi)X,Y)+\eta (\nabla_X\varphi Y)+\eta (\nabla_{\varphi X}Y)-\eta ([\varphi X,Y])=$$
$$=g(\nabla_{\varphi X}\xi,Y)+\eta (\nabla_{\varphi X}Y)-\eta ([\varphi X,Y])=(\pounds_{\varphi X}\eta)Y=(\pounds_{\varphi
Y}\eta)X=$$
$$-g((\pounds_{\xi}\varphi)Y,X)+\eta (\nabla_Y\varphi X)+\eta (\nabla_{\varphi Y}X)-\eta ([\varphi
Y,X]).$$

Hence $g((\pounds_{\xi}\varphi)X,Y)=g((\pounds_{\xi}\varphi)Y,X)$.

For the second statement, using $Proposition~\ref{l1}$, we have
$$2g((\nabla_X\varphi)\xi,Z)=-g(N^{(1)}(\xi,Z),\varphi X)-2d\eta(\varphi Z,X)=-g((\pounds_{\xi}\varphi)Z,X)+$$
$$+2g(Z,X)-2\eta(X)\eta(Z)=-g((\pounds_{\xi}\varphi)X,Z)+2g(Z,X)-2\eta(X)\eta(Z)$$
and hence $\varphi \nabla_X\xi=hX-X+\eta(X)\xi$. Applying
$\varphi$ we obtain
$$\nabla_X\xi=- \varphi X+\varphi hX.$$
To see the anti-commutativity, note that
$$2g(X,\varphi Y)=2d\eta(X,Y)=g(\nabla_X\xi,Y)-g(\nabla_Y\xi,X)=-g(\varphi X,Y)+g(\varphi hX,Y)+$$
$$+g(\varphi Y,X)-g(\varphi hY,X).$$
Therefore $0=g(\varphi hX,Y)+g(Y,h\varphi X)$ giving $h\varphi +
\varphi h=0.$ From the equality $\varphi
\nabla_X\xi=hX-X+\eta(X)\xi$ we get $h\xi=0.$
\end{proof}
\begin{co}
On a paracontact manifold, $\delta \eta =0$, where $\delta$ is the
co-differential.
\end{co}
\begin{lem}\label{l3}
On a paracontact metric manifold we have the formula
\begin{equation}
(\nabla_{\varphi X}\varphi)\varphi
Y-(\nabla_{X}\varphi)Y=2g(X,Y)\xi-(X-hX+\eta(X)\xi)\eta(Y)
\end{equation}
\end{lem}
\begin{proof}
Either using \eqref{f2} 
or by direct differentiation of
$\nabla_Y\xi=- \varphi Y+\varphi hY$, we obtain
\begin{equation}\label{f4}
(\nabla_XF)(\varphi Y,Z)-(\nabla_XF)(Y,\varphi
Z)=\eta(Y)g(X-hX,\varphi Z)+\eta(Z)g(X-hX,\varphi Y).
\end{equation}
Replacing Z by $\varphi Z$ and using \eqref{f2}, 
we get
\begin{equation}\label{f5}
(\nabla_XF)(\varphi Y,\varphi Z)-(\nabla_XF)(Y,Z)=\eta(Y)g(X-hX,
Z)-\eta(Z)g(X-hX,Y)
\end{equation}
Now, since $dF=0$ we have
\begin{equation}\label{f100}
-(\nabla_XF)(Y,Z)-(\nabla_YF)(Z,X)-(\nabla_ZF)(X,Y)-
\end{equation}
$$-(\nabla_XF)(\varphi Y,\varphi Z)-(\nabla_{\varphi Y}F)(\varphi Z,X)-(\nabla_{\varphi Z}F)(X,\varphi Y)+$$
$$+(\nabla_{\varphi X}F)(\varphi Y,Z)+(\nabla_{\varphi Y}F)(Z,\varphi X)+(\nabla_{Z}F)(\varphi X,\varphi Y)+$$
$$+(\nabla_{\varphi X}F)(Y,\varphi Z)+(\nabla_{Y}F)(\varphi Z,\varphi X)+(\nabla_{\varphi Z}F)(\varphi X,Y)=0.$$
Now $(\ref{f4})$, $(\ref{f5})$ and $(\ref{f100})$ give
$$(\nabla_{\varphi X}F)(\varphi Y,Z)-(\nabla_XF)(Y,Z)=-2g(X,Y)\eta(Z)+g(X-hX+\eta(X)\xi,Z)\eta(Y)$$
from which the result follows.
\end{proof}
We recall that
a \emph{paraSasakian manifold} is a normal paracontact metric
manifold.

\begin{thm}\label{t4}
An almost paracontact metric structure $(\varphi,\xi,\eta,g)$ is
paraSasakian if and only if
\begin{equation}\label{f6}
(\nabla_X\varphi)Y=-g(X,Y)\xi+\eta (Y)X
\end{equation}
In particular, a paraSasakian manifold is K-paracontact.
\end{thm}
\begin{proof}
Suppose $(\varphi,\xi,\eta,g)$ is paraSasakian. Then
$Proposition~\ref{l1}$ yields
$$2g((\nabla_X\varphi)Y,Z)=-2d\eta(\varphi Z,X)\eta(Y)+2d\eta(\varphi Y,X)\eta(Z)=2g(-g(X,Y)\xi+\eta(Y)X,Z).$$
Conversely, assume $(\nabla_X\varphi)Y=-g(X,Y)\xi+\eta (Y)X$, set
$Y=\xi$ to get $-\varphi \nabla_X\xi=-\eta(X)\xi+X$. Hence $
\nabla_X\xi=-\varphi X$ and therefore
$$2d\eta(X,Y)=g(\nabla_X\xi,Y)-g(\nabla_Y\xi,X)=2g(X,\varphi Y)$$
showing that $(\varphi,\xi,\eta,g)$ is a paracontact metric
structure.

Now, we calculate
$$N_{\varphi}(X,Y)-2d\eta(X,Y)\xi =[\varphi X,\varphi Y]-\varphi[\varphi X,Y]-\varphi[X,\varphi Y]+\varphi^2[X,Y]-2d\eta (X,Y)\xi=$$
$$=(\nabla_{\varphi X}\varphi)Y-(\nabla_{\varphi Y}\varphi)X-\varphi(\nabla_{X}\varphi)Y+\varphi(\nabla_{Y}\varphi)X-2d\eta (X,Y)\xi=-g(\varphi
X,Y)\xi+$$
$$+\eta(Y)\varphi X+g(\varphi Y,X)\xi-\eta(X)\varphi Y-\eta(Y)\varphi X+\eta(X)\varphi Y-2d\eta (X,Y)\xi=$$
$$=2g(X,\varphi Y)-2d\eta (X,Y)\xi=0.$$
Therefore, $(\varphi,\xi,\eta,g)$ is paraSasakian.
\end{proof}

\section{Curvature of Paracontact manifolds}
In this chapter we discuss some aspects of the curvature of
paracontact manifolds. We begin with some preliminaries concerning
the tensor field h.
\begin{pro}\label{p1}
On a paracontact manifold $M^{2n+1}$ we have the formulas
\begin{gather}\label{f7}
(\nabla_{\xi}h)X=- \varphi X + h^2 \varphi X + \varphi R( \xi , X)
\xi,\\
\label{f8} (R( \xi , X) \xi + \varphi R( \xi ,\varphi X) \xi)=
2\varphi^2 X - 2h^2X
\end{gather}
\end{pro}
\begin{proof}
Using Lemma~\ref{l2},
we calculate
$$R(\xi,X)\xi=\nabla_{\xi}(- \varphi X+\varphi hX)+\varphi [\xi,X]-\varphi h[\xi,X].$$
Applying $\varphi$, taking into account that
$\nabla_{\xi}\varphi=0$, we obtain
$$\varphi R(\xi,X)\xi=-\nabla_X\xi+(\nabla_{\xi}h)X+h\nabla_X\xi.$$
Apply Lemma~\ref{l2} to get \eqref{f7}.
 Multiply \eqref{f7} with $\varphi$ to derive
$$R(\xi,X)\xi=\varphi^2X+\varphi (\nabla_{\xi}h)X-h^2X.$$
Taking into account that
$\varphi R(\xi,\varphi X)\xi=\varphi^2X-\varphi(\nabla_{\xi}h)X-h^2X,$
we get \eqref{f8}.
\end{proof}
\begin{co}
On a paracontact  metric manifold $M^{2n+1}$ the Ricci curvature
in the direction of $\xi$ is given by
\begin{equation}\label{f9}
Ric( \xi , \xi)=-2n + |h|^2
\end{equation}
On a $K$-paracontact metric manifold $M^{2n+1}$ we have $Ric( \xi
, \xi)=-2n$.
\end{co}
\begin{pro}\label{p2}
On a paraSasakian manifold
$$R(X,Y)\xi=\eta(X)Y-\eta(Y)X.$$
\end{pro}
\begin{proof} We calculate
$$R(X,Y)\xi=-\nabla_X\varphi Y+\nabla_Y\varphi X +\varphi
[X,Y]=-(\nabla_X\varphi)Y+(\nabla_Y\varphi)X=$$
$$=\eta(X)Y-\eta(Y)X.$$
\end{proof}
\begin{lem}\label{l4}
The curvature tensor of a paracontact metric manifold satisfies
\begin{gather}\label{f10}
R( \xi ,X,Y,Z)=-(\nabla_XF)(Y,Z)+g(X,(\nabla_Y\varphi
h)Z)-g(X,(\nabla_Z\varphi h)Y),\\\label{f11}
R( \xi ,X,Y,Z)+R( \xi ,X,\varphi Y,\varphi Z)-R( \xi ,\varphi
X,\varphi Y,Z)-R( \xi ,\varphi X,Y,\varphi Z)\\\nonumber
=-2(\nabla_{hX}F)(Y,Z)+2g(X-hX,Z)\eta (Y)-2g(X-hX,Y)\eta (Z).
\end{gather}
\end{lem}
\begin{proof}
Differentiating $\nabla_Z\xi=- \varphi Z+\varphi hZ$, we obtain
$$R(Y,Z)\xi=-(\nabla_Y\varphi)Z+(\nabla_Z\varphi)Y+(\nabla_Y\varphi h)Z-(\nabla_Z\varphi h)Y$$
which, since $dF=0$, yields the first formula \eqref{f10}.
Set
\begin{gather*}A(X,Y,Z)=-(\nabla_XF)(Y,Z)-(\nabla_XF)(\varphi Y,\varphi Z)+(\nabla_{\varphi X}F)(Y,\varphi Z)\\
+(\nabla_{\varphi X}F)(\varphi Y,Z)\\
B(X,Y,Z)=g(X,(\nabla_Y\varphi h)Z)+g(X,(\nabla_{\varphi Y}\varphi h)\varphi Z)-g(\varphi X,(\nabla_{Y}\varphi h)\varphi Z)\\-
g(\varphi X,(\nabla_{\varphi Y}\varphi h)Z).
\end{gather*}
Use \eqref{f10} to see that the left hand side of \eqref{f11} is
equal to
$A(X,Y,Z)+B(X,Y,Z)-B(X,Z,Y)$. The proof of  $Lemma~\ref{l4}$
yields
$$A(X,Y,Z)=-2g(X,Y)\eta(Z)+2g(X,Z)\eta(Y).$$
It is straightforward to show that $\eta((\nabla_{\varphi
Y}h)Z)=g(Y+hY,hZ)$. Rewrite $B$ in the form
$$B(X,Y,Z)=g(X,(\nabla_Y\varphi)hZ)-g(\varphi hX,(\nabla_{\varphi Y}\varphi)Z)-
g(\varphi X,(\nabla_{\varphi Y}\varphi)hZ)-$$
$$-g(\varphi X,h(\nabla_{Y}\varphi)Z)+\eta(X)(\nabla_{\varphi Y}\eta)hZ.$$
Use $Lemma~\ref{l4}$ again to obtain
$$B(X,Y,Z)=-2g(hX,(\nabla_Y\varphi)Z)+2g(hX,Y)\eta(Z)+2g(hY,hZ)\eta(X).$$
Finally, compute $A(X,Y,Z)+B(X,Y,Z)-B(X,Z,Y)$, use $dF=0$ to get
the result.
\end{proof}

Let us fix a local coordinates $(x^1,\dots,x^{2n+1})$. We shall
use the Einstein summation convention.
The equations $(\ref{f82})$, $(\ref{con})$ and
$(\ref{f3})$, in local coordinates, have the expression 
$$\eta_r\xi^r=1,\qquad \varphi_r^i\xi^r=0,\qquad \eta_r\varphi^r_j=0,\qquad
\varphi^i_r\varphi^r_j=\delta^i_j-\xi^i\eta_j,$$
$$g_{rs}\varphi^r_j\varphi^s_k=-g_{jk}+\eta_j\eta_k,\qquad g_{jr}\xi^r=\eta_j.$$
We get using \eqref{f2} that
\begin{gather}
\nabla_i\eta_j-\nabla_j\eta_i=2\varphi_{ij}=2g_{ir}\varphi^r_j,\nonumber\\\label{f12}
\nabla_r\varphi^r_j=2n\eta_j,\qquad
\xi^r\nabla_r\varphi^i_j=0,\quad
\nabla_r\eta_s\varphi_i^r\varphi_j^s=\nabla_j\eta_i,\\\label{f101}
\nabla_r\eta_i\varphi_j^r \quad and \quad
\nabla_i\eta_r\varphi_j^r \quad are \quad symmetric \quad in \quad
i,j.
\end{gather}
Moreover, $Lemma~\ref{l2}$ implies
\begin{equation}\label{f13}
\nabla_i\eta_j=\varphi_{ij}+\varphi_{ir}h^r_j,
\quad
h_{ij}=h_{ji}=g_{jr}h^r_i,\quad \varphi^i_rh^r_j=-h^i_r\varphi^r_j,\qquad h_{ij}\xi^j=0.
\end{equation}
Consequently, $(\ref{f13})$ yields
\begin{equation}\label{f14}
\nabla_r\eta_i\nabla^r\eta_j=-g_{ij}+\eta_i\eta_j-2h_{ij}-h_{ir}h^r_j.
\end{equation}
From the equations $(\ref{f8})$ and $(\ref{f9})$ we also have
\begin{gather}\label{f15}
R_{irsj}\xi^r\xi^s-R_{arsb}\xi^r\xi^s\varphi^a_i\varphi^b_j=-2g_{ij}+2\eta_i\eta_j+2h_{ir}h^r_j,
\\\nonumber Ric(\xi,\xi)=-2n+|h|^2,
\end{gather}
where $|h|^2=g^{ir}g^{js}h_{ij}h_{rs}$, for $h=(h_{ij})$.
\begin{lem}\label{l5}
Let $(M,g,\varphi,\eta,\xi)$ be a paracontact pseudo-Riemannian
manifold. Then the Ricci tensor $Ric$ of the Levi-Chevita
connection satisfies the following relations:
\begin{gather}\label{f16}
Ric_{jr}\xi^r=\nabla_r\nabla_j\xi^r=\nabla_r\nabla^r\eta_j-4n\eta_j,
\\\label{f17}
\varphi^s_j\nabla^r\nabla_r\varphi_{ks}+\varphi^s_k\nabla^r\nabla_r\varphi_{js}=
2\nabla_r\varphi_{sj}\nabla^r\varphi^s_k-Ric_{jr}\xi^r\eta_k-Ric_{kr}\xi^r\eta_j
\\\nonumber+2h_{jr}h^r_k+4h_{jk}+2g_{jk}-2(4n+1)\eta_j\eta_k.
\end{gather}
\end{lem}
\begin{proof}
Contracting
$R^k_{ijr}\xi^r=\nabla_i\nabla_j\xi^k-\nabla_j\nabla_i\xi^k$ with
respect to i and k, we obtain the first equality in $(\ref{f16})$.
To verify the second equality, we observe that
$\nabla^r\nabla_r\eta_j=\nabla^r(2\varphi_{rj})+\nabla^r\nabla_r\eta_r$.
Then use $(\ref{f12})$ to get $(\ref{f16})$. Next, applying the
hyperbolic Laplacian $\nabla^r\nabla_r$ to
$\varphi^s_j\varphi_{ks}=g_{jk}-\eta_j\eta_k$, we obtain
$$\varphi^s_j\nabla^r\nabla_r\varphi_{ks}+
\varphi^s_k\nabla^r\nabla_r\varphi_{js}-2\nabla_r\varphi_{sj}\nabla^r\varphi^s_k=
-\eta_k\nabla^r\nabla_r\eta_j-\eta_j\nabla^r\nabla_r\eta_k-2\nabla_r\eta_j\nabla^r\eta_k.$$
The latter together with $(\ref{f14})$ and $(\ref{f16})$ yields
$(\ref{f17})$.
\end{proof}
The obstruction an almost  paracontact pseudo-Riemannian manifold 
to be a paraSasakian, described in Theorem~\ref{t4}, is the tensor
$P=(P_{rsi})$  defined by 
\begin{equation}\label{f18}
P_{rsi}=\nabla_r\varphi_{si}-\eta_ig_{rs}+\eta_sg_{ri}.
\end{equation}
\begin{lem}\label{l6}On a paracontact metric manifold
$P_{rsi}P^{rs}_j$ is given by
\begin{equation}\label{f19}
P_{rsi}P^{rs}_j
=\nabla_r\varphi_{si}\nabla^r\varphi^s_j+2h_{ij}-g_{ij}-(2n-1)\eta_i\eta_j.
\end{equation}
\end{lem}
\begin{proof}
First we get
$$P_{rsi}P^{rs}_j
=\nabla_r\varphi_{si}\nabla^r\varphi^s_j+\eta_s\nabla_i\varphi^s_j+\eta_s\nabla_j\varphi^s_i+g_{ij}-(2n+1)\eta_i\eta_j.$$
Since $\eta_s\nabla_i\varphi^s_j=-\varphi^s_j\nabla_i\eta_s$,
applying $(\ref{f13})$ to the last equation, we obtain
$(\ref{f19})$.
\end{proof}
We define the *-Ricci tensor $Ric^*_{ij}$ and the *-scalar
curvature $scal^*$ by
$$Ric^*_{ij}=g^{ps}R_{pilk}\varphi^l_j\varphi^k_s, \qquad scal^*=g^{ij}Ric^*_{ij}.$$
\begin{lem}\label{l7}
The symmetric part of the *-Ricci tensor is given by
\begin{equation}\label{f20}
Ric^*_{ij}+Ric^*_{ji}=-Ric_{ij}+Ric_{rs}\varphi^r_i\varphi^s_j-2(2n-1)g_{ij}+
\end{equation}
$$+2(n-1)\eta_i\eta_j+P_{rsi}P^{rs}_j+h_{ir}h^r_j.$$
\end{lem}
\begin{proof}
By the Ricci identity for $\varphi$, we obtain
\begin{equation}\label{f21}
\nabla_l\nabla_k\varphi^i_j-\nabla_k\nabla_l\varphi^i_j=R^i_{lka}\varphi^a_j-R^s_{lkj}\varphi^i_s.
\end{equation}
Contracting the last equation with respect to $i$ and $k$, we get
\begin{equation}\label{f22}
2n\nabla_l\eta_j-\nabla_i\nabla_l\varphi^i_j=-Ric_{la}\varphi^a_j-R^s_{ilj}\varphi^i_s.
\end{equation}
Transvecting $(\ref{f22})$ by $\varphi^l_k$, we obtain
\begin{equation}\label{f23}
2n\nabla_l\eta_j\varphi^l_k-\varphi^l_k\nabla_i\nabla_l\varphi^i_j=-Ric_{la}\varphi^a_j\varphi^l_k+Ric^*_{jk}.
\end{equation}
Transvecting $(\ref{f22})$ by $-\varphi^j_k$, we obtain
\begin{equation}\label{f24}
-2n\nabla_l\eta_j\varphi^j_k+\varphi^j_k\nabla_i\nabla_l\varphi^i_j=Ric_{lk}-Ric_{la}\xi^a\eta_k+Ric^*_{lk}.
\end{equation}
Change $l$ to $j$ in $(\ref{f24})$. Then the obtained result and
$(\ref{f23})$ imply
$$4n\varphi_{rj}\varphi^r_k-\varphi^r_k\nabla^i(\nabla_r\varphi_{ij}-\nabla_j\varphi_{ir})=Ric_{jk}-Ric_{rs}\varphi^r_j\varphi^s_k-Ric_{js}\xi^s\eta_k+2Ric^*_{jk}.$$
Since
$\nabla_r\varphi_{ij}+\nabla_i\varphi_{jr}+\nabla_j\varphi_{ri}=0$,
the above is written as
$$-4n(g_{kj}-\eta_k\eta_j)+\varphi^r_k\nabla^i\nabla_i\varphi_{jr}=Ric_{jk}-Ric_{rs}\varphi^r_j\varphi^s_k-Ric_{js}\xi^s\eta_k+2Ric^*_{jk}.$$
Take the symmetric part of the latter equation, use $(\ref{f17})$
and $(\ref{f19})$ to derive  $(\ref{f20})$.
\end{proof}
We define
$P(X)=(P_{rsi}X^i)$. Then we get
$|P(X)|^2=(P_{rsi}P^{rs}_jX^iX^j)$. By $(\ref{f19})$ it easy to
verify
\begin{equation}\label{f25}
|P(\xi)|^2=|h|^2.
\end{equation}
Therefore, if $(M,\varphi,\eta,g)$ is a K-paracontact manifold,
then $|P(\xi)|^2=0$.

By $Lemma~\ref{l7}$ we obtain the following
\begin{co}
If a paracontact manifold $(M,\varphi,\eta,g)$ is a paraSasakian,
then
\begin{equation}\label{f26}
Ric^*_{ij}+Ric^*_{ji}=-Ric_{ij}+Ric_{rs}\varphi^r_i\varphi^s_j-2(2n-1)g_{ij}+2(n-1)\eta_i\eta_j
\end{equation}
\end{co}
The equalities $(\ref{f20})$ and $(\ref{f9})$ give 
\begin{co}\label{c2}
Let $(M,\varphi,\eta,g)$ be a paracontact manifold. Then
\begin{equation}\label{f27}
scal+scal^* +4n^2= |h|^2 + \frac{1}{2}|\nabla \varphi|^2-2n,
\end{equation}
where $|P|^2=|\nabla \varphi|^2-4n$. If $(M,\varphi,\eta,g)$ is
paraSasakian manifold, then $$scal+scal^* +4n^2=0.$$
\end{co}
In the contact case the identity (\ref{f27}) has been proven by
Olszak (\cite{O1}, see also \cite{T1}).
\begin{thm}\label{t8}
Let $(M,\varphi,\eta,g)$ be a locally conformally equivalent to a
flat paracontact manifold of dimension $2n+1\geqq 5$. For
any unit $X$ orthogonal to $\xi$
\begin{gather}\label{f28}
Ric(X,X)-Ric(\varphi X,\varphi X)=-4n
-\frac{1}{n(2n-3)}(2n(2n+1)+scal)
\\\nonumber+\frac{2n-1}{2n-3}(|P(X)|^2+|h(X)|^2).
\end{gather}
If $(M,\varphi,\eta,g)$ is a conformally flat paraSasakian
manifold and $2n+1\geqq 5$, then
$$Ric(X,X)-Ric(\varphi X,\varphi X)=-4n
-\frac{1}{n(2n-3)}(2n(2n+1)+scal).$$
\end{thm}
\begin{proof}
Recall that a Riemannian manifold is locally conformally flat
exactly when the Weyl curvature vanishes due to the Weyl's
theorem. Let $(M,\varphi,\eta,g)$ be a conformally flat
paracontact manifold. Then the Riemannian curvature tensor $R$ is
expressed as
$$R_{ijkl}=\frac{1}{2n-1}(Ric_{jk}g_{il}-Ric_{ik}g_{jl}-Ric_{jl}g_{ik}+Ric_{il}g_{jk})-
\frac{scal}{2n(2n-1)}(g_{jk}g_{il}-g_{ik}g_{jl})$$
Hence, $Ric^*(X,X)$ for any unit $X \perp \xi$ is given by
$$Ric^*(X,X)=-\frac{1}{2n-1}(Ric(X,X)-Ric(\varphi X,\varphi X))+\frac{scal}{2n(2n-1)}$$
On the other hand, ($\ref{f20}$) gives
$$2Ric^*(X,X)=-Ric(X,X)+Ric(\varphi X,\varphi X))-2(2n-1)+|P(X)|^2+|h(X)|^2.$$
Combining the last two equations we obtain ($\ref{f28}$).
\end{proof}
\begin{rem}
Let $(e_i, \varphi e_i,\xi)$ be an adapted basis of a conformally
flat paracontact manifold. Then, using $(\ref{f9})$ and
$(\ref{f28})$, we can show that the scalar curvature $scal$ is
given by
$$scal=-2n(2n+1)+\frac{2n-1}{4(n-1)(2n-3)}|P|^2+\frac{2n-3}{2(n-1)}|h|^2.$$
\end{rem}
\begin{thm}\label{t9}
If a paracontact manifold $M^{2n+1}$ is of constant sectional
curvature c and dimension $2n+1\geq 5$, then c=-1 and $|h|^2=0$.
\end{thm}
\begin{proof}
Recall from $Proposition~\ref{p1}$ that, $\frac{1}{2}(R( \xi , X)
\xi + \varphi R( \xi ,\varphi X) \xi)=\varphi^2 X - h^2X$; thus if
$R(X,Y)Z=c(g(Y,Z)X-g(X,Z)Y)$, then $c(\eta(X)\xi-X-\varphi^2X)=
2\varphi^2 X -2 h^2X$. Therefore $h^2X=(c+1)\varphi^2X$ and hence
$|h|^2=2n(c+1)$. Now from $Lemma~\ref{l4}$
\begin{gather*}(\nabla_{hX}F)(Y,Z)=-(c+1)g(X,Y)\eta (Z)+(c+1)g(X,Z)\eta (Y)+g(hX,Y)\eta (Z)\\
-g(hX,Z)\eta (Y).
\end{gather*}
Replacing X by $hX$, we have
\begin{gather*}(\nabla_{h^2X}F)(Y,Z)=-(c+1)g(hX,Y)\eta (Z)+(c+1)g(hX,Z)\eta (Y)+g(h^2X,Y)\eta (Z)\\
-g(h^2X,Z)\eta (Y).
\end{gather*}
Hence, $(c+1)((\nabla_{X}\varphi)Y+g(X-hX,Y)\xi-(X-hX)\eta(Y))=0.$ We have two cases

{\bf CASE 1.}
If $c=-1$ then we have $|h|^2=0$.

{\bf CASE 2.}
If $c \not= -1$ then
$(\nabla_{X}\varphi)Y=-g(X-hX,Y)\xi+(X-hX)\eta(Y).$
Using the latter, we compute $|\nabla \varphi|^2$ and applying
$|h|^2=2n(c+1)$, we obtain
$|\nabla \varphi|^2=4n(c+2).$
On the other hand $scal=2n(2n+1)c$ and $scal^*=-2nc$ as is easily
checked. Now from the formula in $Corollary~\ref{c2}$, we obtain
$4n^2(c+1)=4n(c+1).$
This is a contradiction, because $n>1$ and $c\not= -1$.
\end{proof}
We restrict our attention to paraSasakian manifolds. We begin with
\begin{lem}\label{l8}
On a paraSasakian manifold we have
\begin{gather}\label{f29}
R(X,Y,\varphi Z,W)+R(X,Y,Z,\varphi
W)=-d\eta(X,W)g(Y,Z)+d\eta(X,Z)g(Y,W)
\\\nonumber-d\eta(Y,Z)g(X,W)+d\eta(Y,W)g(X,Z),
\\\label{f30}
R(\varphi X,\varphi Y,\varphi Z,\varphi W)-R(X,Y,Z, W) =
\eta(X)\eta(W)g(Y,Z)\\\nonumber+\eta(Y)\eta(Z)g(X,W)
-\eta(Y)\eta(W)g(X,Z)-\eta(X)\eta(Z)g(Y,W),
\\\label{f31}
R(X,\varphi X,Y,\varphi Y)=-R(X,Y,X,Y)+R(X,\varphi Y,X,\varphi
Y)\\\nonumber+\eta(X)\eta(Y)g(X,Y)
+2(d\eta(X,Y)d\eta(X,Y)-g(X,Y)g(X,Y)+|X|^2|Y|^2).
\\\label{f32}
Ric(X,\varphi Y)+Ric(\varphi X,Y)=-d\eta(X,Y)
\end{gather}
\end{lem}
\begin{proof}
The first equality follows by definition and
$(\nabla_X\varphi)Y=-g(X,Y)\xi+\eta (Y)X$. Using the first
equality we obtain the second. The third equality follows by the
first Bianchi identity to $R(X,\varphi X,Y,\varphi Y)$ and using
the first equality. Finally choosing a $\varphi$-basis and second
equality we obtain the fourth equality.
\end{proof}
\begin{co}
On a paraSasakian manifold for X,Y,Z,W orthogonal to $\xi$ we have
$$R(\varphi X,\varphi Y,\varphi Z,\varphi W)=R(X,Y,Z,W);$$
$$R(X,\varphi X,Y,\varphi Y)=-R(X,Y,X,Y)+R(X,\varphi Y,X,\varphi
Y)+$$
$$+2(d\eta(X,Y)d\eta(X,Y)-g(X,Y)g(X,Y)+|X|^2|Y|^2);$$
$$Ric(X,\varphi Y)+Ric(\varphi X,Y)=0.$$
\end{co}
The Bianchi identities and equation \eqref{f30} yield
\begin{lem}\label{l10}
The Ricci tensor $Ric$ of a $(2n+1)-$dimensional paraSasakian
manifold $M$ satisfies the relations 
\begin{gather*}
Ric(X,Y)=\frac{1}{2}\sum_{i=1}^{2n+1}R(X,\varphi Y,e_i,\varphi e_i)-(2n-1)g(X,Y)-\eta(X)\eta(Y),\\
Ric(\varphi X,\varphi Y)=-Ric(X,Y)-2n\eta(X)\eta(Y),\\
(\nabla_ZRic)(X,Y)=(\nabla_XRic)(Y,Z)-(\nabla_{\varphi Y}Ric)(\varphi X,Z)-\eta(X)Ric(\varphi Y,Z)
\\-2\eta(Y)Ric(\varphi X,Z)-2n\eta(X)g(\varphi Y,Z)-4n\eta(Y)g(\varphi X,Z).
\end{gather*}
\end{lem}

\section{Canonical paracontact connection and conformal (gauge) transformation}\label{gau}

Let $(M,\varphi,\xi,\eta,g)$ be a paracontact manifold. All
paracontact forms $\tilde\eta$ generating the same horizontal
distribution $\mathbb D=Ker\,\eta$ are connected by
$\tilde\eta=\sigma\eta$ for a positive smooth function $\sigma$ on
$M$. We consider another paracontact form $\widetilde{\eta}=\sigma
\eta$ and define structure tensors
$(\widetilde{\varphi},\widetilde{\xi},\widetilde{g})$
corresponding to $\widetilde{\eta}$ using the condition:

$(\star)$ For each point $x$ of $M$, the actions of $\varphi$ and
$\widetilde{\varphi}$ are identical  on $\mathbb{D}_x$

By calculating $d\widetilde{\eta}=d(\sigma \eta)$, we obtain
\begin{equation}\label{f55}
2\widetilde{\varphi}_{ij}=\sigma_i\eta_j-\sigma_j\eta_i+2\sigma
\varphi_{ij},
\end{equation}
where $\sigma_i=\nabla_i\sigma$. By
$\widetilde{\xi}^i\widetilde{\varphi}_{ij}=0$,
$\widetilde{\eta}_i\widetilde{\xi}^i=1$ and $(\ref{f55})$, we
obtain $\widetilde{\xi}\sigma=\frac{1}{\sigma}\xi \sigma$, and
\begin{equation}\label{f56}
\widetilde{\xi}^k=\frac{1}{\sigma}\xi^k-\frac{1}{2\sigma^2}\varphi_j^k\sigma^j.
\end{equation}
So we define $\zeta$ by
$\zeta^k=-\frac{1}{2\sigma}\varphi_j^k\sigma^j$ and get
$$\widetilde{\xi}^k=\frac{1}{\sigma}(\xi^k+\zeta^k).$$
By
$\widetilde{\varphi}_{ij}\widetilde{\varphi}^{jk}=\delta_i^k-\widetilde{\xi}^k\widetilde{\eta}_i$
and $\widetilde{\eta}_j\widetilde{\varphi}^{jk}=0$,
$\widetilde{\varphi}^{jk}$ is determined:
\begin{equation}\label{f57}
\widetilde{\varphi}^{jk}=\frac{1}{\sigma}\varphi^{jk}.
\end{equation}
Now, by the condition $(\star)$ we can put
$\widetilde{\varphi}_j^i=\varphi_j^i+v^i\eta_j$ for some vector
field $v$ on $M$. By $\widetilde{\eta}_i\widetilde{\varphi}_j^i=0$
and
$\widetilde{\varphi}_j^i\widetilde{\varphi}_k^j=\delta_k^i-\widetilde{\xi}^i\widetilde{\eta}_k$,
$v$ is determined:
$$v^i=\frac{1}{2\sigma}(\sigma^i-\xi \sigma \cdot \xi^i).$$
By the expressions of $\widetilde{\varphi}_{ij}$ and
$\widetilde{\varphi}_j^i$, we obtain
$$\widetilde{g}_{ij}=\sigma(g_{ij}-\eta_i\zeta_j-\eta_j\zeta_i)+\sigma(\sigma-1+|\zeta|^2)\eta_i\eta_j.$$
The inverse matrix $(\widetilde{g}^{jk})$ of
$(\widetilde{g}_{ij})$ is given by
$$\widetilde{g}^{jk}=\frac{1}{\sigma}(g^{jk}-\xi^j\xi^k)+\frac{1}{\sigma^2}(\xi^j+\zeta^j)(\xi^k+\zeta^k).$$
The last relation can be  rewritten as
\begin{equation}\label{f58}
\sigma(\widetilde{g}^{jk}-\widetilde{\xi}^j\widetilde{\xi}^k)=g^{jk}-\xi^j\xi^k.
\end{equation}
Summarizing the above discussions,  we obtain 
\begin{lem}\label{l11}
Under condition $(\star)$, a gauge transformation $\eta
\rightarrow \widetilde{\eta}=\sigma \eta$ of a paracontact form
$\eta$ induces the transformation of the structure tensors of the
form:
\begin{gather*}
\widetilde{\xi}^k=\frac{1}{\sigma}(\xi^k+\zeta^k),\quad \zeta^k
=-\frac{1}{2\sigma}\varphi_j^k\sigma^j,\\
\widetilde{\varphi}_j^i=\varphi_j^i+\frac{1}{2\sigma}(\sigma^i-\xi \sigma \cdot \xi^i)\eta_j,\\
\widetilde{g}_{ij}=\sigma(g_{ij}-\eta_i\zeta_j-\eta_j\zeta_i)+
\sigma(\sigma-1+|\zeta|^2)\eta_i\eta_j.
\end{gather*}
\end{lem}
We call the transformation of the structure tensors given by
$Lemma~\ref{l11}$ a \emph{gauge (conformal) transformation of
paracontact pseudo-Riemannian structure}. When $\sigma$ is
constant this is a $\mathbb{D}$-homothetic transformation studied
in the $Subsection~\ref{homo}$

On a strongly pseudo-convex CR-manifold Tanaka \cite{Ta1} and
Webster \cite{W} introduced a canonical connection preserving the
structure called \emph{Tanaka-Webster connection}. Tanno
generalized this connection extending its definition to the
general contact metric manifold.

Following \cite{T2}, we consider the connection
$\widetilde{\nabla}$ defined by
\begin{equation}\label{tan-web}
\begin{aligned}
\widetilde{\nabla}_XY=\nabla_XY+\eta(X)\varphi Y-\eta(Y)\nabla_X\xi+(\nabla_X\eta)Y\cdot \xi\\
=\nabla_XY+\eta(X)\varphi Y+\eta(Y)(\varphi X-\varphi hX)+g(X,\varphi Y)\cdot \xi-g(hX,\varphi Y)\cdot \xi.
\end{aligned}
\end{equation}
The torsion of this connection is then
\begin{gather}\label{tprtw}
T(X,Y)=\eta(X)\varphi Y-\eta(Y)\varphi X-\eta(Y)\nabla_X\xi+\eta(X)\nabla_Y\xi+2d\eta(X,Y)\xi\\\nonumber
=\eta(X)\varphi hY-\eta(Y)\varphi hX+2g(X,\varphi Y)\xi.
\end{gather}
\begin{pro}\label{p6}
On a paracontact manifold the  
connection $\widetilde{\nabla}$ 
has the properties
\begin{equation}\label{tnweb}
\begin{aligned}
\widetilde{\nabla}\eta=0,\quad \widetilde{\nabla}\xi=0,\quad \widetilde{\nabla}g=0,\\
(\widetilde{\nabla}_X\varphi)Y=(\nabla_X\varphi)Y+g(X-hX,Y)\xi-\eta(Y)(X-hX),\\
T(\xi,\varphi Y)=-\varphi T(\xi,Y),\quad Y \in \Gamma(\mathbb{D})\quad or \quad Y\in \Gamma(TM)\\
T(X,Y)=2d\eta(X,Y)\xi,\quad X,Y\in \Gamma(\mathbb{D}).
\end{aligned}
\end{equation}
\end{pro}
\begin{proof}
Calculation is straightforward by using \eqref{tan-web} and \eqref{tprtw}.
\end{proof}
\begin{defn}
We call the connection $\widetilde\nabla$ defined above on a paracontact manifold
\emph{the canonical paracontact connection.}
\end{defn}
We calculate the curvature of $\widetilde{\nabla}$. Let $W$ be
the $(1,2)$-tensor field expressing the difference between
$\widetilde{\nabla}$ and $\nabla,\quad W_{ij}^k=\widetilde{\Gamma}_{ij}^k-\Gamma_{ij}^k.$
We obtain using \eqref{tan-web} that
\begin{equation}\label{f59}
\widetilde{R}_{ijk}^l=R_{ijk}^l+\nabla_i\varphi_k^l\eta_j-\nabla_j\varphi_k^l\eta_i+2\varphi_{ij}\varphi_k^l-
\varphi_s^l\nabla_j\xi^s\eta_i\eta_k+\varphi_s^l\nabla_i\xi^s\eta_j\eta_k+
\end{equation}
$$+\xi^l\nabla_i\eta_s\varphi_k^s\eta_j-\xi^l\nabla_j\eta_s\varphi_k^s\eta_i-\xi^lR_{ijk}^s\eta_s-
\eta_kR_{ijs}^l\xi^s+\nabla_j\eta_k\nabla_i\xi^l-\nabla_i\eta_k\nabla_j\xi^l.$$
Contracting \eqref{f59} with respect to $i$ and $l$, we obtain
$$\widetilde{Ric}_{jk}=Ric_{jk} -2g_{jk}+2\eta_j\eta_k-\eta_kRic_{js}\xi^s-R_{jsrk}\xi^s\xi^r-
\nabla_r\eta_k\nabla_j\xi^r.$$
Since $\widetilde{Ric}_{jk}\xi^j\xi^k=0$, we define the scalar
curvature of the canonical paracontact connection of a
paracontact pseudo-Riemannian manifold $(M,\xi,\eta,g)$ by
$W_1=g^{jk}\widetilde{Ric}_{jk}$. Using $(\ref{f16})$, we obtain
$\nabla_r\eta_s\nabla^s\xi^r=-Ric_{rs}\xi^r\xi^s.$ Hence,
\begin{equation}\label{f60}
W_1=scal-Ric(\xi,\xi)-4n.
\end{equation}
Let $f$ and $f'$ be two functions on a paracontact
pseudo-Riemannian manifold $(M,\xi,\eta,g)$. We define operator
$\triangle_{\mathcal{D}}$ acting on the space of functions by
using the hyperbolic Laplacian $\triangle$ and $\xi$:
$$\triangle_{\mathbb{D}}f=\triangle f-\xi \xi f=(g^{ij}-\xi^i\xi^j)\nabla_i\nabla_jf,$$
and $(df;df')_{\mathbb{D}}$ by
$$(df;df')_{\mathbb{D}}=(g^{ij}-\xi^i\xi^j)\nabla_if\nabla_jf'.$$
Furthermore, $|df|_{\mathbb{D}}^2$ means $(df;df)_{\mathbb{D}}$,
which is equal to $|df|^2-( \xi f)^2$.
\begin{thm}
Let $(\eta,g)\rightarrow (\widetilde{\eta}=\sigma
\eta,\widetilde{g})$ be a conformal (gauge) transformation of
paracontact pseudo-Riemannian structure. Then the transformation
of the scalar curvature $W_1$ of the canonical paracontact
connection is given by
\begin{equation}\label{f61}
\sigma
\widetilde{W}_1=W_1-\frac{2(n+1)}{\sigma}\triangle_{\mathbb{D}}\sigma-
\frac{(n+1)(n-2)}{\sigma^2}|d\sigma|_{\mathbb{D}}^2.
\end{equation}
\end{thm}
\begin{proof}
We follow the scheme in \cite{T2}. Geometric quantities
corresponding to $\widetilde{g}$ are denoted by $\sim$. We define
$\widetilde{W}_{jk}^i$ by
$$\widetilde{W}_{jk}^i=\widetilde{\Gamma}_{jk}^i-\Gamma_{jk}^i.$$
Then
$$\widetilde{W}_{jk}^i=\frac{1}{2}\widetilde{g}^{ia}(\nabla_j\widetilde{g}_{ak}+
\nabla_k\widetilde{g}_{aj}-\nabla_a\widetilde{g}_{jk}).$$
The Ricci tensor $\widetilde{Ric}$ is given by
$$\widetilde{Ric}_{jl}=Ric_{jl}+\nabla_r\widetilde{W}_{lj}^r-\nabla_l\widetilde{W}_{rj}^r+
\widetilde{W}_{lj}^s\widetilde{W}_{rs}^r-\widetilde{W}_{rj}^s\widetilde{W}_{ls}^r.$$
Transvecting the last equality by $g^{jl}-\xi^j\xi^l$ and using
$(\ref{f58})$ we obtain
\begin{equation}\label{f62}
\sigma(\widetilde{scal}-\widetilde{Ric}_{jl}\xi^j\xi^l)=
scal-Ric_{jl}\xi^j\xi^l+(g^{jl}-\xi^j\xi^l)\nabla_r\widetilde{W}_{lj}^r\hspace{3cm}
\end{equation}
$$-
(g^{jl}-\xi^j\xi^l)\nabla_l\widetilde{W}_{rj}^r+
(g^{jl}-\xi^j\xi^l)\widetilde{W}_{lj}^s\widetilde{W}_{rs}^r-
(g^{jl}-\xi^j\xi^l)\widetilde{W}_{rj}^s\widetilde{W}_{ls}^r.$$
First we calculate the following:
$$\widetilde{W}_{lj}^r(g^{jl}-\xi^j\xi^l)=\frac{1}{2}\widetilde{g}^{ra}(\nabla_l\widetilde{g}_{aj}+
\nabla_j\widetilde{g}_{al}-\nabla_a\widetilde{g}_{jl})\sigma(\widetilde{g}^{jl}-
\widetilde{\xi}^j\widetilde{\xi}^l)=$$
$$=\sigma\widetilde{g}^{ra}[\nabla_j(\widetilde{\xi}^j\widetilde{\eta}_a)-
\nabla_j(\frac{1}{\sigma}(g^{jl}-\xi^j\xi^l))\widetilde{g}_{al}]-
\frac{1}{2}\sigma\widetilde{g}^{ra}\widetilde{g}_{jl}\nabla(\frac{1}{\sigma}(g^{jl}-\xi^j\xi^l)).$$
After some calculation, we obtain
\begin{gather}\label{f63}
\widetilde{W}_{lj}^r(g^{jl}-\xi^j\xi^l)=\frac{n}{\sigma}\xi \sigma
\cdot \xi^r-\frac{n}{\sigma}\sigma^r,
\\\label{f64}
\nabla_r(\widetilde{W}_{lj}^r(g^{jl}-\xi^j\xi^l))=\frac{n}{\sigma^2}|d\sigma|_{\mathbb{D}}^2-\frac{n}{\sigma}\triangle_{\mathbb{D}}\sigma.
\end{gather}
Next, using 
$2\nabla_r(\xi^j\xi^l)\widetilde{W}_{lj}^r=
\nabla_r\xi^j\xi^l\widetilde{g}^{ra}(\nabla_l\widetilde{g}_{aj}+
\nabla_j\widetilde{g}_{al}-\nabla_a\widetilde{g}_{jl}),$
we derive
\begin{equation}\label{f65}
\nabla_r(\xi^j\xi^l)\widetilde{W}_{lj}^r=4n(\sigma-1+|\zeta|^2)+
\frac{1}{\sigma}(\sigma^r\zeta^s+\sigma^s\zeta^r)\nabla_s\eta_r+
\frac{1}{\sigma^2}\xi\sigma\zeta^r\zeta^s\nabla_s\eta_r+
\end{equation}
$$+\frac{1}{\sigma}\zeta^r\nabla_r\eta_j\zeta^a\nabla^j\eta_a-2\varphi_{rj}\nabla^j\zeta^r+
\frac{1}{\sigma}\zeta^r\nabla_r\eta_j(\nabla^j\zeta_a+\nabla_a\zeta^j)\zeta^a+
\frac{1}{\sigma^2}(1-|\zeta|^2)\sigma^j\zeta^r\nabla_r\eta_j.$$
By a direct calculation we get
\begin{equation}\label{f66}
\widetilde{g}^{ra}\nabla_j\widetilde{g}_{ra}=\frac{2(n+1)}{\sigma}\sigma_j
\end{equation}
Therefore
\begin{equation}\label{f67}
(g^{jl}-\xi^j\xi^l)\nabla_l\widetilde{W}_{rj}^r=\frac{(n+1)}{\sigma}\triangle_{\mathbb{D}}\sigma-\frac{(n+1)}{\sigma^2}|d\sigma|_{\mathbb{D}}^2.
\end{equation}
The fifth term of the right-hand side of $(\ref{f62})$ is
\begin{equation}\label{f68}
(g^{jl}-\xi^j\xi^l)\widetilde{W}_{lj}^s\widetilde{W}_{rs}^r=-\frac{n(n+1)}{\sigma^2}|d\sigma|_{\mathbb{D}}^2
\end{equation}
due to  $(\ref{f63})$ and $(\ref{f66})$.

The sixth term of right-hand side of $(\ref{f62})$ is
$$(g^{jl}-\xi^j\xi^l)\widetilde{W}_{rj}^s\widetilde{W}_{ls}^r=\frac{1}{4}(g^{jl}-\xi^j\xi^l)[-\nabla_j\widetilde{g}_{rs}\nabla_l\widetilde{g}^{rs}-2\nabla_a\widetilde{g}_{jr}\nabla_s\widetilde{g}_{lb}\widetilde{g}^{sa}\widetilde{g}^{rb}+2\nabla_a\widetilde{g}_{jr}\nabla_b\widetilde{g}_{ls}\widetilde{g}^{sa}\widetilde{g}^{rb}].$$
The right-hand side of the last equality is calculated as follows:
$$\frac{1}{4}(g^{jl}-\xi^j\xi^l)\nabla_j\widetilde{g}_{rs}\nabla_l\widetilde{g}^{rs}=-\frac{n+2}{2\sigma^2}|d\sigma|_{\mathbb{D}}^2+\frac{1}{2}(2-\sigma-\frac{1}{\sigma}-|\zeta|^2)|\nabla\xi|^2+(1-\frac{1}{\sigma})(\nabla\xi;\nabla\zeta)+$$
$$+\frac{1}{\sigma}\xi^r\sigma^s\nabla_s\zeta_r-\frac{1}{2\sigma}|\nabla\zeta|^2+\frac{1}{\sigma}\zeta^r\nabla_l\zeta_r\eta_s\nabla^l\zeta^s-\frac{1}{2\sigma}(\sigma+1+|\zeta|^2)\xi^r\nabla_l\zeta_r\eta_s\nabla^l\zeta^s+\frac{1}{2\sigma}|\nabla_{\xi}\zeta|^2,$$
\\
$$\frac{1}{2}(g^{jl}-\xi^j\xi^l)\nabla_a\widetilde{g}_{jr}\nabla_s\widetilde{g}_{lb}\widetilde{g}^{sa}\widetilde{g}^{rb}=\frac{n}{\sigma^2}[|d\sigma|_{\mathbb{D}}^2+(\frac{1}{\sigma}-1)(\xi\sigma)^2]+(\frac{1}{2\sigma}(\sigma-1)^2+|\zeta|^2)|\nabla\xi|^2+$$
$$+\frac{1}{\sigma}(1-\sigma)(\nabla\xi;\nabla\zeta)+(\frac{1}{2\sigma^2}(\sigma-1)^2+\frac{1}{\sigma}|\zeta|^2)|\nabla_{\zeta}\xi|^2+\frac{1}{2\sigma}|\nabla\zeta|^2+\frac{1}{\sigma^2}(1-\sigma)\zeta^a\nabla_a\eta_r(\xi^b+\zeta^b)\nabla_b\zeta^r-$$
$$-\frac{2}{\sigma}\zeta^r\sigma^s\nabla_s\eta_r-\frac{2}{\sigma^2}\xi\sigma\cdot\zeta^r\zeta^s\nabla_s\eta_r+\frac{1}{\sigma}\zeta^s\nabla_a\eta_s\zeta^l\nabla^a\zeta_l+\frac{1}{2\sigma}(1+|\zeta|^2)\zeta^r\nabla_a\eta_r\zeta^s\nabla^a\eta_s+$$
$$+\frac{1}{\sigma^2}\zeta^s\zeta^a\nabla_a\eta_s\zeta^l(\xi^b+\zeta^b)\nabla_b\zeta_l+\frac{1}{2\sigma^2}(1+|\zeta|^2)(\zeta^r\zeta^s\nabla_s\eta_r)^2,$$
\\
$$\frac{1}{2}(g^{jl}-\xi^j\xi^l)\nabla_a\widetilde{g}_{jr}\nabla_b\widetilde{g}_{ls}\widetilde{g}^{sa}\widetilde{g}^{rb}=-\frac{1}{2\sigma^2}|d\sigma|_{\mathbb{D}}^2+\frac{1}{\sigma}\sigma^r\zeta^s\nabla_s\eta_r+\frac{1}{\sigma^2}\xi\sigma(\nabla_r\zeta^r+\zeta^r\zeta^s\nabla_s\eta_r)-$$
$$-\frac{1}{2}|\zeta|^2\nabla_b\eta_a\nabla^a\xi^b+\frac{1}{\sigma}(\nabla_{\zeta}\xi;\nabla_{\xi}\zeta)-\frac{1}{\sigma}(1+|\zeta|^2)\zeta^b\nabla_b\eta_a\zeta^s\nabla^a\eta_s-\frac{1}{\sigma}\zeta^b\nabla_b\eta_a\zeta^r\nabla^a\zeta_r-$$
$$-\frac{1}{2\sigma^2}|\nabla_{\zeta}\xi+\nabla_{\xi}\zeta+\nabla_{\zeta}\zeta|^2-\frac{1}{2\sigma^2}(1+|\zeta|^2)(\zeta^r\zeta^s\nabla_s\eta_r)^2-\frac{1}{\sigma^2}\zeta^r\zeta^a\nabla_a\eta_r\zeta^s(\xi^b+\zeta^b)\nabla_b\zeta_s.$$
Since $\nabla_r\zeta^r=-\frac{n}{\sigma}\xi\sigma$, we obtain
\begin{equation}\label{f72}
(g^{jl}-\xi^j\xi^l)\widetilde{W}_{rj}^s\widetilde{W}_{ls}^r=-\frac{n-1}{2\sigma^2}|d\sigma|_{\mathbb{D}}^2-|\zeta|^2\varphi_{rs}\nabla^s\xi^r+\frac{1}{2}\zeta^r\nabla_l\eta_r\zeta^s\nabla^l\eta_s+
\end{equation}
$$+\frac{1}{\sigma^2}\xi\sigma\cdot\zeta^r\zeta^s\nabla_s\eta_r+\frac{2}{\sigma}\zeta^r\sigma^s\nabla_s\eta_r+\frac{1}{\sigma}(1+|\zeta|^2)\zeta^b\nabla_b\eta_a\zeta^s\nabla^a\eta_s+\frac{1}{\sigma}\zeta^a\nabla_a\eta_r(\nabla_b\zeta^r+\nabla^r\zeta_b)\zeta^b-$$
$$-(\frac{1}{2}-\frac{1}{\sigma}+\frac{1}{\sigma}|\zeta|^2)|\nabla_{\zeta}\xi|^2.$$
Since
$$(\nabla_a\eta_r\nabla_b\xi^r-\nabla_r\eta_a\nabla^r\eta_b)\zeta^a\zeta^b=-2(\varphi_a^r\nabla_b\eta_r+\varphi_a^r\nabla_r\eta_b)\zeta^a\zeta^b=\frac{1}{\sigma}(\nabla_a\eta_b+\nabla_b\eta_a)\zeta^a\zeta^b,$$
etc., summarizing $(\ref{f64})$, $(\ref{f65})$,
$(\ref{f67})$-$(\ref{f72})$, we get
$$\sigma(\widetilde{scal}-\widetilde{Ric}_{jl}\widetilde{\xi}^j\widetilde{\xi}^l)=scal-Ric_{jl}\xi^j\xi^l+4n(\sigma-1)-\frac{2(n+1)}{\sigma}\triangle_{\mathbb{D}}\sigma-\frac{(n+1)(n-2)}{\sigma^2}|d\sigma|_{\mathbb{D}}^2,$$
from which we obtain $(\ref{f61})$.
\end{proof}
\begin{co}\label{c5}
For function $f$ on $M$,
\begin{equation}\label{f69}
\widetilde{\triangle}_{\mathbb{D}}f=\frac{1}{\sigma}\triangle_{\mathbb{D}}f+\frac{n}{\sigma^2}(d\sigma;df)_{\mathbb{D}}
\end{equation}
\end{co}
\begin{proof}
By definition of $\triangle_{\mathbb{D}}$ and
$\widetilde{\triangle}_{\mathbb{D}}$ we obtain
$$\widetilde{\triangle}_{\mathbb{D}}f=(\widetilde{g}^{rs}-\widetilde{\xi}^r\widetilde{\xi}^s)\widetilde{\nabla}_rf_s=\frac{1}{\sigma}(g^{rs}-\xi^r\xi^s)(\nabla_rf_s-\widetilde{W}_{rs}^af_a)=\frac{1}{\sigma}\triangle_{\mathbb{D}}f-\frac{1}{\sigma}(g^{rs}-\xi^r\xi^s)\widetilde{W}_{rs}^af_a.$$
Applying $(\ref{f63})$ to the last line, we get $(\ref{f69})$.
\end{proof}

\subsection{$\mathbb{D}$-homothetic transformations}\label{homo}

In this section we consider homothetic gauge transformation, i.e.
conformal transformation with constant function. Our main
observation here is that these transformations preserve the
$\eta$-Einstein condition in the paraSasakian case. Moreover, we
show that any $\eta$-Einstein paraSasakian manifold with $scal
\not= 2n$ is $\mathbb D$-homothetic to an Einstein manifold.


Set $\sigma=\alpha=const.$ in Lemma~\ref{l11} to get
\begin{equation}\label{f47}
\overline{g}_{jk}=\alpha g_{jk}+\beta \eta_j \eta_k,
\end{equation}
where $\alpha$ and $\beta=\alpha(\alpha-1)$ are constants satisfying $\alpha\not=0$
and $\alpha+\beta>0$. The inverse matrix $(\overline{g}^{ij})$ of
$(\overline{g}_{jk})$ is given by
\begin{equation}\label{f48}
\overline{g}^{ij}=\frac{1}{\alpha}
g^{ij}-\frac{\beta}{\alpha(\alpha+\beta)} \xi^i \xi^j.
\end{equation}
Denoting by $W_{jk}^i$ the difference
$\overline{\Gamma}_{jk}^i-\Gamma_{jk}^i$ of Christoffel symbols,
we have on  a paracontact manifold
\begin{equation}\label{f49}
W_{jk}^i=-\frac{\beta}{\alpha}(\varphi_j^i\eta_k+\varphi_k^i\eta_j)-\frac{\beta}{2(\alpha+\beta)}\xi^i(\nabla_j\eta_k+\nabla_k\eta_j)
\end{equation}
which follows from $(\ref{f47})$ and $(\ref{f48})$.

We assume to the end of this subsection that $M$ is a $K$-paracontact manifold. We have
\begin{equation}\label{f50}
W_{jk}^i=-\frac{\beta}{\alpha}(\varphi_j^i\eta_k+\varphi_k^i\eta_j).
\end{equation}
Substitute \eqref{f50} into
$\overline{R}_{ijk}^l=R_{ijk}^l+\nabla_iW_{jk}^l-\nabla_jW_{ik}^l+W_{is}^lW_{jk}^s-W_{js}^lW_{ik}^s,
$ we obtain
\begin{gather}\label{f52}
\overline{R}_{ijk}^l=R_{ijk}^l-\frac{\beta}{\alpha}(2\varphi_k^l\varphi_{ij}-
\varphi_j^l\varphi_{ik}+\varphi_i^l\varphi_{jk})+
\\\nonumber+\frac{\beta}{\alpha}(\nabla_j\varphi_i^l\eta_k+\nabla_j\varphi_k^l\eta_i-
\nabla_i\varphi_j^l\eta_k-\nabla_i\varphi_k^l\eta_j)+
\frac{\beta^2}{\alpha^2}(\delta_j^l\eta_i\eta_k-\delta_i^l\eta_j\eta_k).
\end{gather}
Contracting with respect $i$ and $l$, we have
\begin{equation}\label{f53}
\overline{Ric}_{jk}=Ric_{jk}+2\frac{\beta}{\alpha}g_{jk}-2\frac{\beta}{\alpha^2}((2n+1)\alpha+n\beta)\eta_j\eta_k,
\end{equation}
where we have used $(\ref{f12})$. Contracting the last equation
with $(\ref{f48})$, we get
\begin{equation}\label{f54}
\overline{scal}=\frac{1}{\alpha}scal+2n\frac{\beta}{\alpha^2}.
\end{equation}
\begin{defn}
If the Ricci tensor of a $K$-paracontact manifold $M$ is of
the form
$$Ric(X,Y)=ag(X,Y)+b\eta(X)\eta(Y),$$
$a$ and $b$ being constant, then $M$ is called an
\emph{$\eta-$Einstein manifold}.
\end{defn}
\begin{pro}\label{t15}
Let $M$ be a $(2n+1)$-dimensional paraSasakian manifold. If the
Ricci tensor $Ric$ of $M$ satisfies
$Ric(X,Y)=ag(X,Y)+b\eta(X)\eta(Y),$ then $a$, $b$  and $scal$ are
constant.
\end{pro}
\begin{proof}
From the assumption on the Ricci tensor $Ric$ we have $a+b=-2n$
and $scal=(2n+1)+b.$ Then we have $Za=-Zb$ and
$Z(scal)=(2n+1)Za+Zb=-2nZb$. On the other hand, $Lemma~\ref{l10}$
implies
$Z(scal)=2Za+2(\xi b)\eta(Z)=-2Zb+(\xi b)\eta(Z).$
Therefore, we obtain
$(n-1)(Zb)=-(\xi b)\eta(Z).$
Put $Z=\xi$ in the latter to find $\xi b=0$. Hance
$Zb=0$, which shows that $b$ is constant. Then $a$ and $scal$ are
also constant.
\end{proof}
\begin{thm}\label{t12}
Let $(M,\varphi,\xi,\eta,g)$ be a paracontact manifold and
$$\overline{\varphi}=\varphi,\quad
\overline{\xi}=\frac{1}{\alpha}\xi, \quad \overline{\eta}=\alpha
\eta,\quad \overline{g}=\alpha g +(\alpha^2-\alpha)\eta \otimes
\eta,\quad \alpha=const.\not=0$$ 
be a $\mathbb D$-homothetic transformation. Then
$(\overline{\varphi},\overline{\xi},\overline{\eta},\overline{g})$
is a paracontact structure too.
\begin{enumerate}
\item[i).] If $(\varphi,\xi,\eta,g)$ is a
$K$-paracontact structure (resp. paraSasakian), then
$(\overline{\varphi},\overline{\xi},\overline{\eta},\overline{g})$
is also a $K$-paracontact structure (resp. paraSasakian).
\item[ii).] If $(\varphi,\xi,\eta,g)$ is
a $\eta$-Einstein paraSasakian structure, then
$(\overline{\varphi},\overline{\xi},\overline{\eta},\overline{g})$
is also a $\eta$-Einstein paraSasakian structure.
\item[iii).] If $(\varphi,\xi,\eta,g)$ is
a $\eta$-Einstein paraSasakian structure with $scal \not= 2n$,
then there exists a constant  $\alpha$ such that
$(\overline{\varphi},\overline{\xi},\overline{\eta},\overline{g})$
is an Einstein paraSasakian structure.
\end{enumerate}
\end{thm}
\begin{proof}
If $\xi$ is a Killing vector field with respect to $g$, then
$\overline{\xi}$ is also a Killing vector field with respect to
$\overline{g}$, since $\xi$ leaves $\eta$ invariant.  The
paraSasakian structure is preserved since the normality conditions is preserved
under the $\mathbb D$-homothetic transformations which proves i).
If $(\varphi,\xi,\eta,g)$ is
$\eta$-Einstein paraSasakian structure then we have
\begin{equation}\label{f79}
Ric(X,Y)=ag(X,Y)+b\eta(X)\eta(Y),
\end{equation}
It follows from  $Proposition~\ref{t15}$  that $a$ and $b$ are
constant and   $(\ref{f9})$ yields
\begin{equation}\label{f80}
a+b=-2n
\end{equation}
Then equality $(\ref{f79})$ has the form
\begin{equation}\label{f81}
Ric(X,Y)=(\frac{scal}{2n}+1)g(X,Y)-(2n+1+\frac{scal}{2n})\eta(X)\eta(Y),
\end{equation}
From  $(\ref{f53})$ and $(\ref{f54})$, for
$\beta=\alpha^2-\alpha$, we derive
\begin{equation}\label{f90}
\overline{Ric}(X,Y)=(\frac{\overline{scal}}{2n}+1)\overline{g}(X,Y)-
(2n+1+\frac{\overline{scal}}{2n})\overline{\eta}(X)\overline{\eta}(Y).
\end{equation}
Finally, we prove iii). If we chose
$\alpha=\frac{2n-scal}{4n^2+4n}$, the equation $(\ref{f54})$ for
$\beta=\alpha^2-\alpha$ gives $\overline{scal}=-2n(2n+1)$. Then we
obtain $\overline{Ric}(X,Y)=-2n\overline{g}(X,Y)$ from
$(\ref{f90})$.
\end{proof}

\subsection{Integrable paracontact manifolds}
Here we consider the case when the paracomplex structure $\varphi$
defined on $\mathbb D$ is formally integrable, i.e. the Nijenhuis
tensor $N_{\varphi}=[\varphi,\varphi]$ satisfies certain
integrability conditions. We see below that in this case the
canonical paracontact connection shares many of the properties of
the Tanaka-Webster connection on $CR$-manifold. We begin with
\begin{defn}
An almost paracontact structure $(\eta,\varphi,\xi)$ is said to be
\emph{integrable} if the almost para-complex structure
$\varphi_{|_{\mathbb D}}$ satisfies the conditions
\begin{equation}\label{new1}
N_{\varphi}(X,Y)=0, \quad X,Y\in \Gamma(\mathbb D).
\end{equation}
and
\begin{equation}\label{mon1}
[\varphi X,Y]+[X,\varphi Y] \in \Gamma(\mathbb D), \quad X,Y\in
\Gamma(\mathbb D).
\end{equation}

Equivalently, the $\pm$-eigendistrubutions $\mathbb D^{\pm}$ of $\varphi$ are
formally integrable in the
sense that
\begin{equation}\label{new2}
[\mathbb D^{\pm},\mathbb D^{\pm}] \in \mathbb D^{\pm}.
\end{equation}
\end{defn}
Indeed, \eqref{new2} means that
\begin{equation}\label{new3}
\begin{aligned}
-\varphi[X+\varphi X,Y+\varphi Y]+[X+\varphi X,Y+\varphi Y]=0, \quad X,Y\in \Gamma(\mathbb D)\\
\varphi[X-\varphi X,Y-\varphi Y]+[X-\varphi X,Y-\varphi Y]=0, \quad X,Y\in \Gamma(\mathbb D)\\
\end{aligned}
\end{equation}
which  is clearly equivalent to \eqref{new1}

%
\begin{thm}\label{t13}
A  paracontact pseudo-Riemannian manifold $(M,g,\varphi,\eta,\xi)$
is integrable if and only if the canonical paracontact connection
preserves the structure tensor
$\varphi$,$$\widetilde{\nabla}\varphi=0$$
\end{thm}
\begin{proof}
It suffices to show that the integrability conditions \eqref{new1}
and \eqref{mon1}
are satisfied if and only if
\begin{equation}\label{new4}(\nabla_X\varphi)Y+g(X-hX,Y)\xi-\eta(Y)(X-hX)=0.
\end{equation}
Indeed, \eqref{new1} can be also written in the form
$$\varphi [X-\eta(X)\xi,\varphi Y]+\varphi [\varphi X,Y-\eta(Y)\xi]=[\varphi X,\varphi Y]+[X-\eta(X)\xi,Y-\eta(Y)\xi],$$
where $X,Y \in \Gamma(TM).$

From the last identity, we obtain
\begin{equation}
g((\nabla_{\varphi X}\varphi) Y,Z)-g((\nabla_{\varphi Y}\varphi)
X,Z)+g((\nabla_{X}\varphi) Y,\varphi
Z)-g((\nabla_{Y}\varphi)X,\varphi Z)-
\end{equation}
$$-2d\eta(X,Y)\xi+\eta(X)(\nabla_Y)Z-\eta(Y)(\nabla_X)Z+\eta(X)g(\nabla_{\varphi
Y}\xi,\varphi Z)-\eta(Y)g(\nabla_{\varphi X}\xi,\varphi Z)-$$
$$-\eta(X)(g(\nabla_{\xi}\varphi Y,\varphi Z)+g(\nabla_{\xi}Y,Z))+\eta(Y)(g(\nabla_{\xi}\varphi X,\varphi Z)+g(\nabla_{\xi}X,Z))=0$$
From the identity $\nabla_{\xi}\varphi=0$ and the
Lemma$~\ref{l2}$, we get
$$g((\nabla_{\varphi X}\varphi) Y,Z)-g((\nabla_{\varphi Y}\varphi)
X,Z)+g((\nabla_{X}\varphi) Y,\varphi
Z)-g((\nabla_{Y}\varphi)X,\varphi Z)+$$
$$+\eta(X)((\nabla_Y)Z+(\nabla_Z)Y)-\eta(Y)((\nabla_X)Z+(\nabla_Z)X)-2d\eta(X,Y)\xi=0$$
That is,
\begin{equation}\label{f74}
\varphi_k^h\nabla_i\varphi_{hj}-\varphi_k^h\nabla_j\varphi_{hi}+\varphi_i^s\nabla_s\varphi_{kj}-\varphi_j^s\nabla_s\varphi_{si}+\eta_i(\nabla_j\eta_k+\nabla_k\eta_j)-
\end{equation}
$$-\eta_j(\nabla_i\eta_k+\nabla_k\eta_i)-2\varphi_{ij}\eta_k=0.$$
Since $d\eta$ is closed, the third and the fourth terms of the
left-hand side of $(\ref{f74})$ are calculate as follows:
$$\varphi_i^s\nabla_s\varphi_{kj}-\varphi_j^s\nabla_s\varphi_{si}=-\varphi_i^s(\nabla_k\varphi_{js}+\nabla_j\varphi_{sk})+\varphi_j^s(\nabla_k\varphi_{is}+\nabla_i\varphi_{sk})=\nabla_i(\eta_j\eta_k)-$$
$$-\nabla_j(\eta_i\eta_k)+\nabla_k(\eta_i\eta_j)-2\varphi_j^s\nabla_k\varphi_{si}+\varphi_k^s\nabla_j\varphi_{si}-\varphi_k^s\nabla_i\varphi_{sj}.$$
Therefore $(\ref{f74})$ is equivalent to
\begin{equation}\label{f75}
\varphi_j^s\nabla_k\varphi_{si}-\eta_i\nabla_k\eta_j=0.
\end{equation}
From the identity $(\ref{f75})$ it follows that
$$\widetilde{\nabla}_i\varphi_{kj}=\nabla_i\varphi_{kj}+
\nabla_i\eta_s\varphi_j^s \eta_k-\nabla_i\eta_s \varphi_k^s \eta_j=
\nabla_i\varphi_{kj}+\varphi_j^s\varphi_s^r\nabla_i\varphi_{rk}-\varphi_k^s\nabla_i\eta_s\eta_j$$
$$=\nabla_i\varphi_{kj}+\nabla_i\varphi_{jk}-\eta_j(\xi^r\nabla_i\varphi_{rk}+
\varphi_k^s\nabla_i\eta_s)
=-\eta_j\nabla_i(\eta_s\varphi_k^s)=0.$$
\end{proof}
The obstruction an integrable paracontact manifold to be normal,
i.e. paraSasakian, is encoded in the (horizontal) torsion
$T(\xi,X)=h(X), X\in \mathbb D$ of the canonical paracontact
connection. The main result here is the following
\begin{thm}\label{parsas}
The torsion of the canonical paracontact connection vanishes on an
integrable paracontact manifold if and only if it is paraSasakian.
\end{thm}
\begin{proof}
Let be $(M,g,\eta)$ an integrable paracontact manifold. From the
$Theorem~\ref{t13}$ follows $\widetilde{\nabla}\varphi=0$. Since
equation $(\ref{tprtw})$ we get
\begin{equation}\label{f83}
T_{ijk}=-\eta_i\varphi_{jk}+\eta_j\varphi_{ik}-\eta_j\nabla_i\eta_k+\eta_i\nabla_j\eta_k+2\varphi_{ij}\eta_k.
\end{equation}
From the last equation we calculate
\begin{equation}\label{f84}
\xi^iT_{ijk}=-\varphi_{jk}+\nabla_j\eta_k=\frac{1}{2}(\nabla_j\eta_k+\nabla_k\eta_j).
\end{equation}
If $M$ is a paraSasakian manifold, then
$P_{rsi}=\nabla_r\varphi_{si}-\eta_ig_{rs}+\eta_sg_{ri}=0$. From
the last equation we calculate
\begin{equation}\label{f85}
\varphi_s^i\nabla_r\eta_i=g_{rs}-\eta_r\eta_s.
\end{equation}
Transvecting $(\ref{f85})$ by $\varphi_k^s$, we obtain
\begin{equation}\label{f86}
\nabla_r\eta_k=\varphi_{rk}.
\end{equation}
From the equation $(\ref{f86})$ we get
$$\xi^iT_{ijk}=\frac{1}{2}(\nabla_j\eta_k+\nabla_k\eta_j)=\frac{1}{2}(\varphi_{rk}+\varphi_{kr})=0.$$
If $\xi^iT_{ijk}=0$ from equation $(\ref{f84})$ we have
\begin{equation}\label{f87}
\nabla_j\eta_k=\varphi_{jk}.
\end{equation}
From equations $(\ref{f75})$ and $(\ref{f87})$ we calculate
$$P_{rsi}=\nabla_r\varphi_{si}-\eta_ig_{rs}+\eta_sg_{ri}=\xi^l\eta_s\nabla_r\varphi_{li}+\eta_i\varphi_s^k\nabla_r\eta_k-\eta_ig_{rs}+\eta_sg_{ri}=-\xi^l\varphi_i^l\nabla_r\eta_l+$$
$$+\eta_i\varphi_s^k\nabla_r\eta_k-\eta_ig_{rs}+\eta_sg_{ri}=-\eta_s\varphi_{rl}\varphi_i^l+\eta_i\varphi_{rk}\varphi_s^k-\eta_ig_{rs}+\eta_sg_{ri}=\eta_s(-g_{ri}+\eta_r\eta_i)-$$
$$-\eta_i(-g_{rs}+\eta_r\eta_s)-\eta_ig_{rs}+\eta_sg_{ri}=\eta_ig_{rs}-\eta_sg_{ri}-\eta_ig_{rs}+\eta_sg_{ri}=0.$$
Therefore $P=0$, which equivalent to $M$ to be a paraSasakian
manifold.
\end{proof}

\section{Paracontact manifolds with torsion}

If we introduce the forms
\begin{equation}\label{f33}
dF^-(X,Y,Z)=dF(\varphi X,Y,Z)+dF(X,\varphi Y,\varphi Z)+
\end{equation}
$$+dF(\varphi
X,Y,\varphi Z)+dF(X,Y,Z);$$
\begin{equation}\label{f34}
dF^{\varphi}(X,Y,Z)=-dF(\varphi X,\varphi Y,\varphi Z).
\end{equation}
and a direct consequence of the definitions and $Proposition~\ref{l1}$
is the following
\begin{pro}\label{p3}
On any almost paracontact manifold the identities hold:
\begin{equation}\label{f35}
dF^-(X,Y,Z)=-N^{(1)}(X,Y,\varphi Z)-N^{(1)}(Y,Z,\varphi
X)-N^{(1)}(Z,X,\varphi Y);
\end{equation}
\begin{equation}\label{f36}
N^{(1)}(X,Y,Z)=N^{(1)}(\varphi X,\varphi
Y,Z)+\eta(Y)N^{(1)}(X,\xi,Z)+\eta(X)N^{(1)}(\xi,Y,Z);
\end{equation}
\begin{equation}\label{f37}
N^{(1)}(X,Y)=(\nabla_{\varphi X}\varphi)Y-(\nabla_{\varphi
Y}\varphi)X+(\nabla_{X}\varphi)\varphi Y -
(\nabla_{Y}\varphi)\varphi X-
\end{equation}
$$-\eta(X)\nabla_Y\xi+\eta(Y)\nabla_X\xi.$$
\end{pro}
\begin{defn}
A linear connection $\overline{\nabla}$ is said to be an almost
paracontact connection if it preserves the almost paracontact
structure:
$$\overline{\nabla} g=\overline{\nabla} \eta=\overline{\nabla} \varphi=0.$$
\end{defn}

\begin{thm}\label{t10}
Let $(M^{(2n+1)},\varphi,\xi,\eta,g)$ be an almost paracontact
metric manifold. The following conditions are equivalent:

1) The tensor $N^{(1)}$ is skew-symmetric and $\xi$ is a Killing
vector field.

2) There exists an almost paracontact linear connection
$\overline{\nabla}$ with totaly skew-symmetric torsion tensor T.

Moreover, this connection is unique and determined by
$$g(\overline{\nabla}_XY,Z)=g(\nabla_XY,Z)+\frac{1}{2}T(X,Y,Z),$$
where the torsion T is defined by
$$T=2\eta \wedge d\eta+d^{\varphi}F-N^{(1)}+\eta \wedge (\xi \lrcorner N^{(1)}).$$
\end{thm}
\begin{proof}
Let assume that such a connection exists. Then
$$0=g(\nabla_X\xi,Z)+\frac{1}{2}T(X,\xi,Z)$$
holds and the skew-symmetric of T yields that $\xi$ is a Killing
vector field, $2d\eta=\xi \lrcorner T$, $\xi \lrcorner d\eta =0$
and
$$T(\varphi X,\varphi Y,Z)+T(\varphi X,Y,\varphi Z)+T(X,\varphi Y,\varphi Z)+T(X,Y,Z)=-N^{(1)}(X,Y,Z).$$
The letter formula shows that $N^{(1)}$ is skew-symmetric. Since
$\varphi$ is $\overline{\nabla}$-parallel, we can express the
Riemannian covariant derivative of $\varphi$ by the torsion form:
$$T( X,\varphi Y,Z)+T( X,Y,\varphi Z)=-2g((\nabla_X\varphi)Y,Z).$$
Taking the cyclic sum in the above equality, we obtain
$$\sigma_{X,Y,Z}T(X,Y,\varphi Z)=-\sigma_{X,Y,Z}g((\nabla_X\varphi)Y,Z)=dF(X,Y,Z).$$
Adding this result to the formula expressing the tensor $N^{(1)}$
by the torsion $T$, come calculations yield
$$T(\varphi X,\varphi Y,\varphi Z)=-dF(X,Y,Z)-g(N^{(1)}(X,Y),\varphi Z)+\eta(Z)N^{(2)}(X,Y).$$
By replacing X,Y,Z by $\varphi X,\varphi Y,\varphi Z$ and using
the symmetry property of the tensor $N^{(1)}$ in
$Proposition~\ref{p3}$, we obtain the formula for the torsion
tensor $T$.

For the converse, suppose that the almost paracontact structure
has properties $1)$ and define the connection $\overline{\nabla}$
by the formulas $2)$. Clearly $T$ is skew-symmetric and
$2d\eta=\xi \lrcorner T=2\nabla \eta$. Since $\xi$ is a Killing
vector field, we conclude $\nabla g=\nabla \xi =0$. Furthermore,
using the conditions $1)$ and $Proposition~\ref{p3}$, we obtain
$\xi \lrcorner dF=N^{(2)}$. Finally we have to prove that $\nabla
\varphi =0$. This follows by straightforward computations using
the relation between $\nabla \varphi$ and the torsion tensor $T$,
$Proposition~\ref{p3}$, as well as the following lemma.
\end{proof}
\begin{lem}\label{l9}
Let $(M^{(2n+1)},\varphi,\xi,\eta,g)$ be an almost paracontact
metric manifold with a totally skew-symmetric tensor $N^{(1)}$.
Then the following equalities hold:
\begin{equation}\label{f38}
\nabla_{\xi}\xi=\xi \lrcorner d\eta=0;
\end{equation}
\begin{equation}\label{f39}
(\nabla_{X}\eta)Y+(\nabla_{Y}\eta)X=-(\nabla_{\varphi
X}\eta)\varphi Y-(\nabla_{\varphi Y}\eta)\varphi X;
\end{equation}
\begin{equation}\label{f40}
N^{(1)}(\varphi X,Y,\xi)=N^{(1)}(X,\varphi Y,\xi)=-N^{(2)}(X,Y)=
\end{equation}
$$=dF(X,Y,\xi)=dF(\varphi X,\varphi Y,\xi).$$
\end{lem}
\begin{proof}
The identities follow from $Proposition~\ref{l1}$, $Lemma~\ref{l3}$ and
formula $(\ref{f36})$.
\end{proof}
We discuss these results for some special paracontact structures.
\begin{thm}\label{t11}
Let $(M^{(2n+1)},\varphi,\xi,\eta,g)$ be an almost paracontact
metric manifold with totally skew-symmetric tensor $N^{(1)}=0$.
Then the condition $dF=0$ implies $N^{(1)}=0$.

1) A paracontact metric structure $(F=d\eta)$ admits an almost
paracontact connection with totally skew-symmetric torsion if and
only if it is paraSasakian. In this case, the connection is
unique, its torsion is given by
$$T=2\eta \wedge d\eta$$
and T is parallel, $\overline{\nabla}T=0$.

2)  A normal $(N^{(1)}=0)$ paracontact metric structure admits a
unique almost paracontact connection with totally skew-symmetric
torsion if and only if $\xi$ is Killing vector field. The torsion
T is given by
$$T=2\eta \wedge d\eta+d^{\varphi}F.$$
\end{thm}
\begin{proof}
If $dF=0$, $Lemma~\ref{l9}$ implies that $N^{(2)}=\xi \lrcorner
N^{(1)}=0$. Then Proposition leads to
$0=dF^-(X,Y,Z)=-3N^{(1)}(\varphi X,Y,Z)$. The assertion that
$\overline{\nabla}T=0$ in a paraSasakian manifold follows by
direct verification.
\end{proof}
We introduce the forms
\begin{equation}\label{f41}
\rho^{\overline{\nabla}}(X,Y)=\frac{1}{2}\sum_{i=1}^{2n+1}R^{\overline{\nabla}}(X,Y,e_i,\varphi
e_i);
\end{equation}
\begin{equation}\label{f42}
t(X)=\frac{1}{2}\sum_{i=1}^{2n+1}T(X,e_i,\varphi e_i);
\end{equation}
\begin{equation}\label{f43}
dt(X,Y)=\frac{1}{2}\sum_{i=1}^{2n+1}dT(X,Y,e_i,\varphi e_i);
\end{equation}
\begin{pro}\label{p4}
Let $(M^{(2n+1)},\varphi,\xi,\eta,g)$ be an almost paracontact
metric manifold with totally skew-symmetric tensor $N^{(1)}$ and
Killing vector $\xi$. Let $\overline{\nabla}$ be the unique almost
paracontact connection with totally skew-symmetric torsion. Then
one has
\begin{equation}\label{f44}
\rho^{\overline{\nabla}}(X,Y)=Ric^{\overline{\nabla}}(X,\varphi
Y)+(\overline{\nabla}_Xt)Y+\frac{1}{2}dt(X,Y)
\end{equation}
\end{pro}
\begin{proof}
We follow the scheme in \cite{F1} and use the curvature properties
of $\overline{\nabla}$ in to calculate $dt(X,Y)$:
$$dt(X,Y)=(\overline{\nabla}_Xt)Y-(\overline{\nabla}_Yt)X+\sigma^T(X,Y,e_i,\varphi e_i)-(\overline{\nabla}_{\varphi e_i}T)(X,Y,e_i).$$
The first Bianchi identity for $\overline{\nabla}$ together with
the latter identity implies
$$4\rho^{\overline{\nabla}}(X,Y)+2Ric^{\overline{\nabla}}(Y,\varphi
X)-2Ric^{\overline{\nabla}}(X,\varphi
Y)=2dt(X,Y)+2(\overline{\nabla}_Xt)Y-2(\overline{\nabla}_Yt)X.$$

Using the relation between the curvature tensors of $\nabla$ and
$\overline{\nabla}$, we obtain
$$Ric^{\overline{\nabla}}(Y,\varphi
X)+Ric^{\overline{\nabla}}(X,\varphi
Y)=-(\overline{\nabla}_Xt)Y-(\overline{\nabla}_Yt)X.$$

The last two equalities lead to the desired formula.
\end{proof}
\begin{pro}\label{p5}
Let $(M^{(2n+1)},\varphi,\xi,\eta,g)$ be a paraSasakian metric
manifold and $\nabla$ be the unique almost paracontact connection
with totally skew-symmetric torsion. Then one has
\begin{equation}\label{f45}
\rho^{\overline{\nabla}}(X,\varphi
Y)=Ric^{\overline{\nabla}}(X,Y)+4(n-1)(g(X,Y)-\eta(X)\eta(Y))
\end{equation}
Moreover, the 2-form $\rho^{\overline{\nabla}}=0$ if and only if
\begin{equation}\label{f46}
Ric(X,Y)=-2(2n-1)g(X,Y)+2(n-1)\eta(X)\eta(Y).
\end{equation}
\end{pro}
\begin{proof}
On a paraSasakian manifold $T=2\eta \wedge d\eta=2\eta \wedge F$
and $\nabla T=0$, where $F(X,Y)=g(X,\varphi Y)$ is the fundamental
form of the paraSasakian structure. Consequently, we calculate
that
$$\overline{\nabla} t=0,dt=8(n-1)F,\sum_{i=1}^{2n+1}g(T(X,e_i),T(Y,e_i))=-8g(X,Y)-8(n-1)\eta(X)\eta(Y).$$

Using the relation between the curvature tensors of $\nabla$ and
$\overline{\nabla}$, we obtain
$$Ric(X,Y)=Ric^{\overline{\nabla}}(X,Y)-2g(X,Y)-2(n-1)\eta(X)\eta(Y)$$
and
$$\rho(X,\varphi Y)=\rho^{\overline{\nabla}}(X,\varphi Y)-(2n-1)(g(X,Y)-\eta(X)\eta(Y)),$$
and the proof follows from $Proposition~\ref{p4}$.
\end{proof}

\textbf{Acknowledgement} This work was partially supported by
Contract 062/2007 with the University of Sofia "St. Kl. Ohridski".

\bibliographystyle{hamsplain}






\end{document}